\documentclass[12pt,a4paper]{amsart}

\usepackage[pdftex,pagebackref]{hyperref}
\hypersetup{letterpaper=true,colorlinks=true}
\hypersetup{pdfpagemode=UseNone,urlcolor=blue}
\hypersetup{linkcolor=blue,citecolor=blue,pdfstartview=FitH}

\usepackage{amsmath,amsfonts,amssymb,amsthm}
\usepackage{graphicx}
\usepackage{enumerate}
\usepackage{verbatim}

\usepackage{color}

\usepackage{ulem}

\def\emph#1{\textit{#1}}

\newenvironment{acknowledgements}{\textbf{Acknowledgements:}}{}

\renewcommand{\geq}{\geqslant}
\renewcommand{\leq}{\leqslant}

\newcommand{\ra}{\rightarrow}
\newcommand{\lra}{\longrightarrow}

\newcommand{\hdot}{{\:\raisebox{3pt}{\text{\circle*{1.5}}}}}

\def\C{{\mathbb C}}
\renewcommand{\P}{\mathbb{P}}

\newcommand{\cO}{\mathcal{O}}
\newcommand{\cE}{\mathcal{E}}
\newcommand{\cF}{\mathcal{F}}
\newcommand{\cI}{\mathcal{I}}
\newcommand{\cZ}{\mathcal{Z}}
\newcommand{\cW}{\mathcal{W}}
\newcommand{\cX}{\mathcal{X}}
\newcommand{\cA}{\mathcal{A}}
\newcommand{\cB}{\mathcal{B}}
\newcommand{\cM}{\mathcal{M}}

\newtheorem{theorem}{Theorem}[section]
\newtheorem{lemma}[theorem]{Lemma}

\newtheorem{corollary}[theorem]{Corollary}
\newtheorem{proposition}[theorem]{Proposition}

\newtheorem{conjecture}[theorem]{Conjecture}

\renewcommand{\hom}{\textrm{Hom}}

\newcommand{\ext}{\textrm{Ext}}

\newcommand{\rk}{\textrm{rank}}

\newcommand{\aut}{\textrm{Aut}}
\newcommand{\cext}{\mathcal{E}\mathit{xt}_{\cX/\P^1}}
\newcommand{\chom}{\mathcal{H}\mathit{om}_{\cX/\P^1}}
\newcommand{\rhom}{\textrm{RHom}}
\newcommand{\crhom}{\mathcal{R}\mathcal{H}\mathit{om}_{\cX/\P^1}}
\renewcommand{\sec}{\mathcal{S}\!\mathit{ec}}
\newcommand{\hor}{\mathcal{H}\!\mathit{or}}

\title{A Note on Formality and Singularities of Moduli Spaces}
\date{\today}
\author{Ziyu Zhang}
\address{Institute for Mathematics, University of Mainz, Staudingerweg 9, 55099 Mainz, Germany}
\email{zhangzy@uni-mainz.de}

\begin{document}

\maketitle

\begin{abstract}
This paper studies formality of the differential graded algebra $\rhom^\hdot(E,E)$, where $E$ is a semistable sheaf on a K3 surface. The main tool is Kaledin's theorem on formality in families. For a large class of sheaves $E$, this DG algebra is formal, therefore we have an explicit description of the singularity type of the moduli space of semistable sheaves at the point represented by $E$. This paper also explains why Kaledin's theorem fails to apply in the remaining case. 
\end{abstract}

\tableofcontents

\section{Introduction}

Moduli spaces of semistable sheaves on complex projective K3 surfaces have been studied intensively for a long time. They are constructed by Geometric Invariant Theory, and share many great properties with the base surfaces. It is known in \cite[Corollary 0.2]{Muk-symp} that the locus parameterizing stable sheaves in such a moduli space is smooth and carries a natural non-degenerate holomorphic 2-form. In particular, when the Mukai vector is primitive, there are no strictly semistable sheaves, therefore the moduli space is a holomorphic symplectic manifold. However, the picture is not so complete when the Mukai vector is non-primitive. 

We require some deformation theory of sheaves in this discussion. (For example, section 3.1 of \cite{KL-local} gives a brief account in the context that we need here.) To get a basic idea, let $(X,H)$ be a polarized complex projective K3 surface and $v$ is a fixed Mukai vector. We consider the moduli space $M_{X,H}(v)$ of semistable sheaves with such a Mukai vector. The closed points of the moduli space $M_{X,H}(v)$ are one-to-one correspondent to polystable sheaves, which are direct sums of stable summands. At the points where the corresponding sheaf $E$ have at least two direct summands, $M_{X,H}(v)$ is singular. The singularity type is related to the solutions of the Maurer-Cartan equation associated to the differential graded algebra (``DG algebra" for short) $\rhom^\hdot(E,E)$. More precisely, we look at a small open neighborhood $U_0$ of $0$ in the complex vector space $\ext^1(E,E)$. The DG algebra $\rhom^\hdot(E,E)$ determines a closed subvariety $U_E$ in $U_0$, which is called a versal deformation space of $E$. The automorphism group $\aut (E)$ acts on the complex vector space $\ext^1(E,E)$, preserving the versal deformation space $U_E$. In fact, since the action of the scalars is trivial, we actually have an action on the vector space $\ext^1(E,E)$ by $G=\aut(E)/\C^*$. Then locally near the point represented by the polystable sheaf $E$, the moduli space is exactly the affine algebraic variety obtained as the quotient $U_E/\!/G$.

 In general it is hard to write down the precise equations for the versal deformation space $U_E$. However, when the DG algebra $\rhom^\hdot(E,E)$ is formal, namely, quasiisomorphic to its homology algebra $\ext^\hdot(E,E)$, by a theorem of Goldman and Millson \cite[Theorem 5.3]{GM-htpy}, the singularities appearing in the versal deformation space $U_E$ has a particularly simple form. In this case, the Yoneda pairing on the complex vector space $\ext^1(E,E)$ gives a quadratic equation $$Q: \ext^1(E,E)\lra \ext^2(E,E)_0,$$ where $\ext^2(E,E)_0$ is the trace-free part of the vector space $\ext^2(E,E)$. The zero locus $Q^{-1}(0)\cap U_0$ is a versal deformation space $U_E$ of $E$, which has at worst quadratic singularities. 

Furthermore, it's also known, for example in \cite[Section 2.7]{KLS}, that the quadratic map $Q$ is actually the symplectic moment map of the $G$ action on $V$. So locally near the point represented by $E$, the moduli space looks exactly like the quotient of the zero fiber of the moment map by the group action. Therefore both the versal deformation space $U_E$ and the local moduli space $U_E/\!/G$ can be described in a very explicit way, which is critical in a complete classification of singularities appearing in all moduli spaces of semistable sheaves on K3 surfaces. 

The following conjecture was essentially raised in \cite{KL-local}:

\begin{conjecture}\label{mainconj}
Let $X$ be a complex projective K3 surface and $H$ be a generic polarization. 
Let $E$ be a polystable sheaf. Then the DG algebra $\rhom^\hdot(E,E)$ is formal. 
\end{conjecture}

If this conjecture was proven true, then we can give a complete list of singularities types (which turns out to be a countable list) which appear in the moduli spaces of semistable sheaves on K3 surfaces. 

In \cite{KL-local}, the following special case was proven by Kaledin and Lehn: 

\begin{proposition}
Conjecture \ref{mainconj} holds when $E=\cI_Z^{\oplus n}$ where $\cI_Z$ is the ideal sheaf of a 0-dimensional subscheme $Z$. 
\end{proposition}

The proof of the above special case made essential use of a theorem of Kaledin \cite[Theorem 4.3]{Kal-family} about formality in families. 

In this paper, we will try to extend the above proof and apply Kaledin's theorem in a more general setting. To fix the notation, let $(X,H)$ be a generically polarized K3 surface as above. Let 
\begin{equation*}\label{eqn_st_decomp}
E=\bigoplus_{i=1}^m E_i^{\oplus n_i}
\end{equation*}
be the stable decomposition of the polystable sheaf $E$, where $E_i$'s are distinct stable summands. The genericness of the choice of the polarization $H$ guarantees that the Mukai vectors of all direct summands $E_i$ are proportional to that of $E$. Then we have the following property: 

\begin{proposition}\label{mainprop}
Let $(X,H)$ be a generically polarized K3 surface as above. Let $v=(r,c,a)$ be a fixed Mukai vector such that $r\geq 1$. We further assume there's at least one $\mu$-stable sheaf on $X$ with Mukai vector $v$. Let $E$ be a polystable sheaf with Mukai vector $v$ and $$E=\bigoplus_{i=1}^m E_i^{\oplus n_i}$$ be the stable decomposition of the $E$, where $E_i$'s are distinct stable summands. Let $r_i$ be the rank of $E_i$. Then Conjecture \ref{mainconj} holds in the following cases: 
\begin{enumerate}
\item when $r_i\geq 2$ for all $i$;
\item when $r_i=1$ for all $i$.
\end{enumerate}
\end{proposition}

As a simple corollory of proposition \ref{mainprop}, we have

\begin{corollary}
Assume all the assumptions of proposition \ref{mainprop}. 
Let $G=\textrm{Aut}(F)/\mathbb{C}^*$,
$V=\textrm{Ext}^1(F,F)$, and $W=\textrm{Ext}^2(F,F)_0$. Let $\gamma:V\rightarrow W$ be the Yoneda pairing which is equivariant under the conjugation actions of $G$. Then in either of the two cases in proposition \ref{mainprop}, there exists
an analytic open neighborhood of $[F]$ in $M_{X,H}(v)$ which is isomorphic
to an analytic open neighborhood of $0$ in the GIT quotient
$\gamma^{-1}(0)/\!/G$.
\end{corollary}

We will also see that, if some of the stable summands are of rank 1 and others are of rank at least 2, then there's a new phenomenon popping up, and Kaledin's theorem cannot be applied directly. 

The paper is organized as follows. In section 2 we will recall Kaledin's key theorem on formality in families. In section 3 we will deal with the first case, which primarily concerns the hyperholomorphic bundles. The second case will be proved in section 4, which mainly deals with ideal sheaves of points. In section 5 we will see the new phenomenon in the mixed cases, which prevents us from applying Kaledin's theorem. 

\begin{acknowledgements}
I express here my deep gratitudes to my doctoral advisor Jun Li for his consistent help and encouragement in the past few years. I would like to thank Daniel Huybrechts and Manfred Lehn for many helpful discussions and suggestions, and for providing me excellent working conditions in SFB-TR 45 program. I would also like to thank Eyal Markman, Dmitry Kaledin, Misha Verbitsky and K\={o}ta Yoshioka for kindly answering my questions and pointing out important references. Moreover, thank Xiaowei Wang for his initial suggestion to me of looking into this problem. This work is supported by SFB-TR 45. 
\end{acknowledgements}

\section{Kaledin's Theorem}\label{sec_kal}

Kaledin's theorem on formality in families is the following: 

\begin{theorem}\cite[Theorem 4.3]{Kal-family}\label{kalthm}
Let $\cA^\hdot$ be a DG algebra of flat quasi-coherent sheaves on a reduced irreducible scheme $X$. Let $\cB^\hdot$ be the homology algebra of the DG algebra $\cA^\hdot$. Assume that the sheaves $\cB^\hdot$ are coherent and flat on $X$, and that for any integers $l$, $i$, the degree-$l$ component $\mathcal{HH}^i_l(\cB^\hdot)$ of the $i$-th Hochschild cohomology sheaf $\mathcal{HH}^i(B)$ is also coherent and flat. 
\begin{enumerate}[(i)]
\item Assume that $X$ is affine. If the fiber $\cA_x^\hdot$ is formal for a generic point $x\in X$, then it is formal for an arbitrary point $x\in X$;
\item Assume that $\mathcal{HH}^2_l(\cB^\hdot)$ has no global sections for all $l\leq -1$. Then the DG algebra $\cA_x^\hdot$ is formal for every point $x\in X$. \qed
\end{enumerate}
\end{theorem}

Our goal is to prove the DG algebra $A(E)=\rhom^\hdot(E,E)$ is formal. Its homology algebra is given by $B(E)=\ext^\hdot(E,E)$. We decompose the polystable sheaf $E$ into $$E=\bigoplus_{i=1}^m E_i^{\oplus n_i}$$ where $E_i$'s are distinct stable sheaves. We denote the Mukai vector of $E_i$ by $v_i$ and the moduli space of stable sheaves of Mukai vector $v_i$ by $\mathcal{M}(v_i)$. If we deform a certain $E_i$ in the moduli space $\mathcal{M}(v_i)$ over a smooth connected base $S$, we will get a flat deformation $\cE$ of $E$ over $S$. A resolution of $\cE$ over $X\times S$ will result in a locally constant deformation of the DG algebra $\rhom^\hdot(E,E)$ and its homology algebra $\ext^\hdot(E,E)$ over $S$. Therefore the Hochschild cohomology sheaves in the assumption of Theorem \ref{kalthm} (i) are also locally constant. Thus Theorem \ref{kalthm} (i) applies. Since for every $i$, the moduli space $\mathcal{M}(v_i)$ is smooth connected, we conclude that, to prove Conjecture \ref{mainconj}, it suffices to prove it for the case that each $E_i$ is a generic closed point in the corresponding moduli space $\mathcal{M}(v_i)$. 

Now we want to decide what a generic sheaf in the moduli space $\mathcal{M}(v_i)$ looks like, under the assumption of Proposition \ref{mainprop}. We assumed that there's at least one $\mu$-stable sheaf with Mukai vector $v$. The following criterion guarantees that in each of $\cM(v_i)$, there's also a $\mu$-stable sheaf. 

\begin{lemma}\cite[Remark 2.2]{Yos-irr}\cite[Lemma 4.4, Remark 4.3]{Yos-examples}
Let $(X,H)$ be a generic polarized K3 surface. Let $v=(lr, l\xi, a)$ be a Mukai vector, such that $r$ and $\xi$ are coprime. Then there's at least one $\mu$-stable sheaf with Mukai vector $v$ unless the following two conditions simultaneously hold: 
\begin{itemize}
\item $\displaystyle \frac{\xi^2+2}{2r}$ is an integer;
\item $\langle v^2 \rangle <2l^2$. \qed
\end{itemize}
\end{lemma}

It is easy to see that, if $v$ is multiplied by a scalar, both of the conditions remain the same. We also know that the Mukai vectors $v_i$ of the stable summands $E_i$ in the decomposition of $E$ is rationally proportional to $v$. Therefore, under the assumption in Proposition \ref{mainprop} that there exists a $\mu$-stable sheaf with Mukai vector $v$, we conclude that the moduli space $\cM(v_i)$ contains at least one $\mu$-stable sheaf as well. Take into consideration that $\mu$-stability is an open condition, we know that there is an open subset in $\cM(v_i)$ which parameterizes $\mu$-stable sheaves. Moreover, the following lemma of Yoshioka says that a generic sheaf in the $\mu$-stable locus is locally free. 

\begin{lemma}\cite[Remark 3.2]{Yos-abelian}\label{st+lf}
Let $r$ be the rank component of the Mukai vector $v$. Then the complement of the locus of locally free sheaves in the moduli of $\mu$-stable sheaves $\cM^{\mu\textrm{-st}}(v)$ has codimension $r-1$. \qed
\end{lemma}

Put the above lemmas together, we can conclude the following: 

\begin{corollary}\label{genericity}
Under the assumption of Proposition \ref{mainprop}, for every Mukai vector $v_i$ of the summand $E_i$ in the decomposition of $E$, the moduli space $\cM(v_i)$ of stable sheaves with Mukai vector $v_i$ has an open subset, which parameterizes 
\begin{enumerate}
\item $\mu$-stable locally free sheaves for $r_i\geq 2$;
\item sheaves of the form $L\otimes \cI_Z$ for $r_i=1$, where $L$ is a line bundle and $\cI_Z$ is an ideal sheaf of 0-dimensional subscheme $Z$ of $X$. 
\end{enumerate}
\end{corollary}

\begin{proof}
When $r_i\geq 2$, in the above argument we have seen that, $\cM(v_i)$ has an open subset which parameterizes $\mu$-stable sheaves, due to the fact that $v_i$ is proportional to $v$ and $\mu$-stability is an open condition. Furthermore, since local freeness is also an open condition, by Lemma \ref{st+lf} we can find an open subset in $\mu$-stable locus which parameterizes $\mu$-stable locally free sheaves. 

When $r_i=1$, by the fact that every sheaf on a smooth surface can be embedded in its double dual with a quotient sheaf with 0-dimensional support, we conclude that the moduli space $\cM(v_i)$ itself is the open set that we are looking for.
\end{proof}

\section{Case of Hyperholomorphic Bundles}\label{sec_hyper}

In this section we will prove Proposition \ref{mainprop} (1). By the argument in the previous section, to prove it for any $E$ of a certain decomposition type, it suffices to prove it for the case that each $E_i$ is generic point in the corresponding moduli space $\cM(v_i)$. By Corollary \ref{genericity}, without loss of generality we can assume every $E_i$ is a $\mu$-stable locally free sheaf. 

To prove Proposition \ref{mainprop} (1), we recall the beautiful idea in \cite{KL-local} of the twistor space of a K3 surface. Roughly speaking, starting from a K3 surface $X$ together with a hyperk\"{a}hler metric, we can produce a complex manifold $\cX$, which is diffeomorphic to $X\times\P^1$ but not holomorphic. The projection $\pi:\cX\ra\P^1$ is holomorphic and each fiber of $\pi$ is equipped with a complex structure which is naturally induced by the hyperk\"{a}hler metric on $X$. In particular, $\pi^{-1}(0)$ is the K3 surface $X$ we start with. 

We also need the notion of a hyperholomorphic bundle \cite[Definition 2.1]{Ver-hyper}. Assume we have a bundle with a Hermitian connection over a K3 surface $X$ with a hyperk\"{a}hler metric. The connection is called hyperholomorphic if it is integrable with respect to the complex structure in any fiber of the twistor family associated to the hyperk\"{a}hler structure on $X$. 

A similar notion on K\"{a}hler manifolds is the Hermitian-Yang-Mills connection \cite{UY-HYM}. Let $B$ be a holomorphic bundle on a K\"{a}hler manifold $M$ with a holomorphic Hermitian connection $\theta$ and a curvature $\Theta\in\Lambda^{1,1}\otimes\textrm{End}(B)$. Then the connection is called Hermitian-Yang-Mills if $$\Lambda(\Theta)=\lambda\cdot\textrm{Id}$$ where $\Lambda$ is a Hodge operator and $c$ is a constant. Then we have the following theorem: 

\begin{theorem}\cite{UY-HYM}\label{UY}
A stable holomorphic vector bundle over a compact K\"{a}hler manifold admits a unique Hermitian-Yang-Mills connection. \qed
\end{theorem}

Now we are ready to prove the following lemma: 

\begin{lemma}\label{dualhyper}
Let $X$ be a K3 surface with a chosen hyperk\"{a}hler metric. Let $E$ be a polystable sheaf on $X$ whose stable summands are all locally free. Then $E^\vee\otimes E$ is a hyperholomorphic bundle. 
\end{lemma}

\begin{proof}
Let $$E=\bigoplus_{i=1}^m E_i^{\oplus n_i}$$ be the stable decomposition as above. By Theorem \ref{UY}, every $E_i$ admits a Hermitian metric $h_i$ and a Hermitian-Yang-Mills connection $\theta_i$ with its curvature form $\Theta_i$, such that $$\Lambda(\Theta_i)=\lambda_i\cdot\textrm{Id}.$$ Since the first Chern class $$c_1(E_i)=\frac{\sqrt{-1}}{2\pi}\textrm{Tr}(\Theta_i),$$ we have $$c_1(E_i)\wedge\omega=\lambda_i\cdot\rk(E_i)\cdot\omega\wedge\omega$$ where $\omega$ is the K\"{a}hler form. By the assumption that the polarization is generic, we know that $c_1(E_i)$ and $\rk(E_i)$ are proportional for all $i$. Therefore all the $\lambda_i$'s are equal and denoted by $\lambda$. This implies that the natural Hermitian metric $h$ and Hermitian connection $\theta$ with the curvature form $\Theta$ is Hermitian-Yang-Mills. It again naturally induces a Hermitian connection $\theta^\vee$ with its curvature form $\Theta^\vee$ on $E^\vee$, where $$\Theta^\vee=-\Theta^\textrm{T}.$$ Furthermore, $E^\vee\otimes E$ is also naturally equipped with a Hermitian metric, a Hermitian connection whose curvature form is given by $$\Omega=\textrm{Id}\otimes\Theta+\Theta^\vee\otimes\textrm{Id}.$$ Since we have $$\Lambda(\Theta)=\lambda\cdot\textrm{Id},$$ we know that $$\Lambda(\Theta^\vee)=-\lambda\cdot\textrm{Id}$$ and $$\Lambda(\Omega)=0.$$ By the same argument as in the proof of \cite[Theorem 2.4]{Ver-hyper}, we know that $\Omega$ is of type $(1,1)$ with respect to any complex structure induced by the hyperk\"{a}hler structure on $X$. By Newlander-Nirenberg Theorem \cite[Proposition 4.17]{Kob-diff}, $E^\vee\otimes E$ is a hyperholomorphic bundle.
\end{proof}

The next tool we need is the twistor transform introduced in \cite{KV-nonHYM}. Given a twistor space $\pi:\cX\ra\P^1$, we can equip it with a natural non-holomorphic projection to the fiber direction $\sigma:\cX\ra X$. Given a bundle $B$ with a connection $\theta$, we can lift it to the total space of the twistor family to obtain a vector bundle with connection $(\sigma^*B, \sigma^*\theta)$. Take the $(0,1)$ part of the connection and we obtain a functor of twistor transform: $$\sigma^*:(B,\theta)\lra(\sigma^*B,(\sigma^*\theta)^{1,0}).$$ The following theorem was proved in \cite{KV-nonHYM}:

\begin{theorem}\cite[Theorem 5.12]{KV-nonHYM}\label{twistlift}
The functor $\sigma^*$ gives an equivalence of the following categories
\begin{enumerate}[(i)]
\item Hyperholomorphic vector bundles on $X$;
\item Holomorphic vector bundles $H$ on $\cX$ such that for all $x\in X$, the restriction of $H$ to a holomorphic horizontal section $$\sigma^{-1}(x)\cong\P^1\subset\cX$$ is a trivial vector bundle $\cO_{\P^1}^{\oplus n}$.
\end{enumerate}\qed
\end{theorem}

In particular, it shows that any hyperholomorphic bundle can be lifted to the twistor space as a holomorphic vector bundle. It turns out that the cohomology groups of the bundle on the fibers of the twistor space has the following amazing property: 

\begin{proposition}\cite[Proposition 6.3]{Ver-coherent}\label{weight}
Let $\pi: \cX\ra\P^1$ be a twistor space of a K3 surface $X$, and $B$ a hyperholomorphic bundle on $X$. Let $\sigma^*B$ be the twistor transform of $B$, a holomorphic vector bundle on $\cX$. Then $$R^i\pi_*(\sigma^*B)\cong\cO_{\P^1}(i)\otimes_{\C}H^i(X,B).$$\qed
\end{proposition}

Now we come back to prove Proposition \ref{mainprop} (1). We follow closely the proof of Proposition 3.1 in \cite{KL-local}. As in \cite[Definition 6.2]{KL-local}, we call a coherent sheaf on $\P^1$ is of weight $l$ if it is a direct sum of several copies of $\cO_{\P^1}(l)$. 

\begin{proof}[Proof of Proposition \ref{mainprop} (1)]
By Lemma \ref{dualhyper}, we know that $F=E^\vee\otimes E$ is a hyperholomorphic bundle. Then by Theorem \ref{twistlift}, we can lift $F$ to the twistor family $\cX$, denoted by $\cF$. Now we consider the sheaf of DG algebras $\crhom^\hdot(\cO_{\cX},\cF)$ on the base of the twistor family $\P^1$. The cohomology sheaves $\cB^\hdot=\cext^\hdot(\cO_{\cX}, \cF)$ is locally constant as a sheaf of graded algebras, so that the Hochschild cohomology sheaves $\mathcal{HH}^\hdot(\cB^\hdot)$ are locally trivial. Therefore Theorem \ref{kalthm} applies. 

Note that for every $i$, $\cB^i=R^i\pi_*\cF$. By Proposition \ref{weight}, $\cB^i$ is of weight $i$. Therefore the degree-$l$ component of the Bar-resolution $\chom^\hdot(\cB^{\hdot\otimes k},\cB^\hdot)$ is of weight $l$ for any integer $l$ and any $k\geq 0$. Since the Hochschild cohomology sheaves $\mathcal{HH}^\hdot_l(\cB^\hdot)$ are computed as the cohomology sheaves of the complex formed by the weight $l$ components of $\chom^\hdot(\cB^{\hdot\otimes k},\cB^\hdot)$, they are themselves of weight $l$. In particular, for any $l\leq -1$, $\mathcal{HH}^2_l(\cB^\hdot)$ has no global sections. We apply Theorem \ref{kalthm} (ii) to conclude that the DG algebra on the fiber of the twistor space $\rhom^\hdot(\cO_X,F)$ is formal. Or equivalently, $\rhom^\hdot(E,E)$ is formal. Hence Proposition \ref{mainprop} (1) is true.
\end{proof}

\section{Case of Ideal Sheaves}

In this section, we will prove Proposition \ref{mainprop} (2). By the assumption, in the stable decomposition of $E$, every summand $E_i$ is a rank 1 sheaf. Since all the Mukai vectors $v_i$'s are proportional, we see that they are all actually the same. Without loss of generality, we can assume that every factor $E_i$ is the ideal sheaf of a 0-dimensional subscheme of length $p$ (possibly 0). 

In fact, fix an $i_0$ and let $L=E_{i_0}^\vee$. Then it's easy to see that for every $E_i$, $E_i\otimes L$ is a rank 1 sheaf with trivial first Chern class on a K3 surface, therefore is the ideal sheaf of a 0-dimensional subscheme. Moreover, the formality of $\rhom^\hdot(E,E)$ is equivalent to the formality of $\rhom^\hdot(E\otimes L, E\otimes L)$. Hence we can safely replace $E$ by $E\otimes L$. 

We can further assume that each $E_i$ is the ideal sheaf of a 0-dimensional subscheme of some length $p\geq 1$. Otherwise, $E_i$ must be the trivial bundle $\cO_X$ and $E=\cO_X^{\oplus n}$ for some positive integer $n$, therefore is a hyperholomorphic bundle. By exactly the same argument as in section \ref{sec_hyper}, $\rhom^\hdot(\cO_X^{\oplus n}, \cO_X^{\oplus n})$ is formal. So we are only left with the case for the length $p\geq 1$. 

We have one more reduction to make. By the same argument as in Section \ref{sec_kal} and applying Theorem \ref{kalthm} (i), we can assume that all 0-dimensional subschemes arising from summands of $E$ are unions of $p$ distinct closed points. 

Before moving on, we want to slightly abuse the notion of the twistor transform defined in \cite{KV-nonHYM} to allow ideal sheaves of distinct points, as in \cite{KL-local}. Let $Z=\cup_{i=1}^p\{z_i\}\subset X$ be the union of $p$ distinct points. For each $i$, $z_i$ determines a horizontal section $\cZ_i=\sigma^*\{z_i\}\cong\P^1$, where $\sigma:\cX\ra X$ is the canonical non-homomorphic projection. Let $\cZ=\cup_{i=1}^p\cZ_i$ be the union of the horizontal sections. We call the ideal sheaf of $\cZ$ in the total space of the twistor family $\cI_{\cZ}$ the twistor transform of $\cI_Z$. Note that on each fiber of the twistor family, it is simply the ideal sheaf of $Z$. 

We want to prove the formality of the DG algebra $\rhom^\hdot(E,E)$. As we have assumed that each stable summand $E_i$ of $E$ is the ideal sheaf of $p$ distinct points, we can lift $E_i$ to the twistor space, denoted by $\cE_i$. Then $$\cE=\bigoplus_{i=1}^m \cE_i^{\oplus n_i}$$ is the twistor transform of $E$. If we consider the corresponding sheaf of DG algebras $\crhom^\hdot(\cE,\cE)$ on $\P^1$, as in the previous case, we still have the local triviality of the cohomology sheaves $\cext^\hdot(\cE,\cE)$ and the Hochschild cohomology sheaves $\mathcal{HH}^m_l(\cB^\hdot)$ for any $m$ and $l$. Now the problem is to analyze the weight of $\cext^k(\cE,\cE)$, for $k=0,1,2$. 

\begin{lemma}\label{case02}
$\cext^k(\cE,\cE)$ is of weight $k$ for $k=0,2$. 
\end{lemma}

\begin{proof}
For $k=0$, due to the stability of each $E_i$, obviously we have $$\cext^0(\cE_i,\cE_i)=\cO_{\P^1}.$$ Therefore it's clear that 
\begin{eqnarray*}
\cext^0(\cE,\cE) &=& \bigoplus_{i=1}^m\cext^0(\cE_i,\cE_i)\otimes\textrm{End}(n_i)\\
&=& \bigoplus_{i=1}^m\cO_{\P^1}\otimes \textrm{End}(n_i).
\end{eqnarray*}
For $k=2$, we need to first analyze the relative canonical bundle $K_{\cX/\P^1}$. It's easy to see that $K_{\cX/\P^1}$ is trivial on each fiber, thus is obtained by pulling back a line bundle from the base $\P^1$. Note that any horizontal section of the twistor family has a normal bundle $\cO_{\P^1}(1)^{\oplus 2}$, therefore a conormal bundle $\cO_{\P^1}(-1)^{\oplus 2}$. Observe that the restriction of $K_{\cX/\P^1}$ on such a horizontal section is exactly the determinant line bundle of the conormal bundle. Therefore we have $$K_{\cX/\P^1}=\pi^*\cO_{\P^1}(-2).$$
Applying relative Serre duality, for each $i$ we have 
\begin{eqnarray*}
\cext^2(\cE_i,\cE_i) &=& \cext^0(\cE_i,\cE_i\otimes K_{\cX / \P^1})^\vee\\
&=& \cext^0(\cE_i,\cE_i\otimes \pi^*\cO_{\P^1}(-2))^\vee\\
&=& (\cext^0(\cE_i,\cE_i)\otimes \cO_{\P^1}(-2)))^\vee\\
&=& \cO_{\P^1}(2).
\end{eqnarray*}
Hence 
\begin{eqnarray*}
\cext^2(\cE,\cE) &=& \bigoplus_{i=1}^m\cext^2(\cE_i,\cE_i)\otimes\textrm{End}(n_i)\\
&=& \bigoplus_{i=1}^m\cO_{\P^1}(2)\otimes \textrm{End}(n_i).
\end{eqnarray*}
\end{proof}

From now on we focus on the middle cohomology sheaf $\cext^1(\cE,\cE)$. It can be decomposed into $$\cext^1(\cE,\cE) =\bigoplus_{i,j=1}^m\cext^1(\cE_i,\cE_j) \otimes\textrm{End}(n_i,n_j).$$ Depending on whether $i$ and $j$ are equal or not, there are 2 types of direct summands which could appear on the right hand side, which are: 
\begin{enumerate}[(i)]
\item $\cext^1(\cI_{\cZ},\cI_{\cZ})$ where $\cI_{\cZ}$ is the twistor transform of an ideal sheaf of $p$ distinct points $\cI_Z$;
\item $\cext^1(\cI_{\cZ},\cI_{\cW})$ where $\cI_{\cZ}$ and $\cI_{\cW}$ are twistor transforms of ideal sheaves of two disjoint set of $p$ points, $\cI_Z$ and $\cI_W$. 
\end{enumerate}

Here is the result in the first case: 

\begin{proposition}\label{samez}
In case (i), $\cext^1(\cI_{\cZ},\cI_{\cZ})$ has weight 1.
\end{proposition}

\begin{proof}
This is part of \cite[Lemma 6.3]{KL-local}.
\end{proof}

The second case is much more complicated. In fact, we don't know a complete answer in this case. The following result is the best that we can prove. Before stating the next proposition, we need one more notation. Let $X^p_0$ be the open subscheme of the product of $p$ copies of $X$, which parameterizes $p$ distinct closed points in $X$. Then in the second case we have $Z,W\in X^p_0$. 

\begin{proposition}\label{diffz}
There exists an analytic open subset $V$ of $X^p_0\times X^p_0$, such that for any pair $(Z,W)\in V$, the sheaf $\cext^1(\cI_{\cZ},\cI_{\cW})$ is of weight 1. More precisely, $\cext^1(\cI_{\cZ},\cI_{\cW})=\cO_{\P^1}(1)^{\oplus(2p-2)}$. 
\end{proposition}

We will prove this proposition in three steps. The first step is to show that, for any $Z$ and $W$ in $X^p_0$, the sheaf $\cext^1(\cI_{\cZ},\cI_{\cW})$ is always an extension of $\cO_{\P^1}(2)^{\oplus (p-1)}$ by $\cO_{\P^1}^{\oplus (p-1)}$. In the second step, we will show the relative extension sheaf $\cext^1(\cI_{\cZ},\cI_{\cW})$ is indeed of weight 1 in a special case. In the third step, we will use a version of upper semi-continuity theorem to finish the proof. 

The first step goes as follows: 

\begin{lemma}\label{lemma_ext}
For any $Z$ and $W$ in $X^p_0$, $\cext^1(\cI_{\cZ},\cI_{\cW})$ is an extension of $\cO_{\P^1}(2)^{\oplus(p-1)}$ by $\cO_{\P^1}^{\oplus(p-1)}$. 
\end{lemma}

\begin{proof}
We resolve $\cI_{\cW}$ by the structure sequence of $\cW$, that is
$$ 0 \lra \cI_{\cW} \lra \cO_{\cX} \lra \cO_{\cW} \lra 0.$$
Then we have the long exact sequence
\begin{eqnarray*}
0 &\lra& \chom(\cI_{\cZ},\cI_{\cW}) \lra \chom(\cI_{\cZ},\cO_{\cX}) \lra \chom(\cI_{\cZ},\cO_{\cW})\\
&\lra& \cext^1(\cI_{\cZ},\cI_{\cW}) \lra \cext^1(\cI_{\cZ},\cO_{\cX}) \lra \cext^1(\cI_{\cZ},\cO_{\cW})\\
&\lra& \cext^2(\cI_{\cZ},\cI_{\cW}) \lra \cext^2(\cI_{\cZ},\cO_{\cX}) \lra \cext^2(\cI_{\cZ},\cO_{\cW}) \lra 0
\end{eqnarray*}

By fiberwise analysis we can easily see that $$\chom(\cI_{\cZ},\cI_{\cW})=\cext^2(\cI_{\cZ},\cI_{\cW})=0.$$
Moreover, we have 
$$\chom(\cI_{\cZ},\cO_{\cX})=\chom(\cO_{\cX},\cO_{\cX})=\cO_{\P^1}$$
as well as
$$\chom(\cI_{\cZ},\cO_{\cW})=\chom(\cO_{\cX},\cO_{\cW})=\cO_{\P^1}^{\oplus p},$$
and the map between them is given by the evaluation of the global functions on $\cW$. Therefore we see that the image of $\chom(\cI_{\cZ},\cO_{\cW})$ in $\cext^1(\cI_{\cZ},\cI_{\cW})$ is $\cO_{\P^1}^{\oplus(p-1)}$. 
For the term $\cext^1(\cI_{\cZ},\cO_{\cW})$, we resolve the ideal sheaf $\cI_{\cZ}$ by 
$$ 0 \lra \cI_{\cZ} \lra \cO_{\cX} \lra \cO_{\cZ} \lra 0$$
and look at the corresponding long exact sequence, part of which is 
$$\cext^1(\cO_{\cX},\cO_{\cW}) \lra \cext^1(\cI_{\cZ},\cO_{\cW}) \lra \cext^2(\cO_{\cZ},\cO_{\cW}).$$
Obviously, $$\cext^1(\cO_{\cX},\cO_{\cW})=0=\cext^2(\cO_{\cZ},\cO_{\cW}),$$ 
therefore $$\cext^1(\cI_{\cZ},\cO_{\cW})=0.$$
Finally, for the term $\cext^1(\cI_{\cZ},\cO_{\cX})$, we use the above structure sequence of $\cZ$ to get a long exact sequence, part of which is
$$
\cext^1(\cO_{\cX},\cO_{\cX}) \lra \cext^1(\cI_{\cZ},\cO_{\cX})
\lra \cext^2(\cO_{\cZ},\cO_{\cX}) \lra \cext^2(\cO_{\cX},\cO_{\cX}).
$$
By fiberwise analysis we easily see $$\cext^1(\cO_{\cX},\cO_{\cX})=0.$$ By Serre duality we have $$\cext^2(\cO_{\cZ},\cO_{\cX}) =(\chom(\cO_{\cX},\cO_{\cZ})\otimes\omega_{\cX/\P^1})^\vee=\cO_{\P^1}(2)^{\oplus p}$$
and $$\cext^2(\cO_{\cX},\cO_{\cX}) =(\chom(\cO_{\cX},\cO_{\cX})\otimes\omega_{\cX/\P^1})^\vee=\cO_{\P^1}(2).$$
Furthermore, since the natural evaluation map $$\chom(\cO_{\cX},\cO_{\cX}) \lra \chom(\cO_{\cX},\cO_{\cZ})$$ is injective, still by Poincare duality, the map $$\cext^2(\cO_{\cZ},\cO_{\cX}) \lra \cext^2(\cO_{\cX},\cO_{\cX})$$ is surjective, whose kernel, by the above analysis, is $\cO_{\P^1}^{\oplus(p-1)}$. 

In summary, we have the short exact sequence $$0 \lra \cO_{\P^1}^{\oplus (p-1)} \lra \cext^1(\cI_{\cZ},\cI_{\cW}) \lra \cO_{\P^1}(2)^{\oplus (p-1)} \lra 0.$$

\end{proof}

By the above lemma, we can conclude that $$\cext^1(\cI_{\cZ},\cI_{\cW})=\cO_{\P^1}^{\oplus a}\oplus \cO_{\P^1}(1)^{\oplus b}\oplus\cO_{\P^1}(2)^{\oplus a}$$ where $a+b=p-1$. Now the goal is to show that $a=0$. In the second step, we show that it is indeed true in a very special case. 

\begin{lemma}\label{manysame}
Let $Z,W\in X^p_0$, such that their intersection is a 0-dimensional subscheme $M$ of length $p-1$. Let $a$ and $b$ be the points in $Z$ and $W$ respectively which are not in $M$. Let $\cI_{\cZ}$ and $\cI_{\cW}$ be the twistor transform of corresponding ideal sheaves. Then $$\cext^1(\cI_{\cZ},\cI_{\cW})=\cO_{\P^1}(1)^{\oplus (2p-2)}.$$
\end{lemma}

\begin{proof}

Let $\cI_{\cM}$ be the twistor transform of the ideal sheaf of $M$. Let $\cA$ and $\cB$ be the horizontal section in the twistor space $\cX$ corresponding to the points $a$ and $b$ respectively. From the short exact sequence $$0 \lra \cI_{\cZ} \lra \cI_{\cM} \lra \cO_{\cA} \lra 0$$ we have
\begin{eqnarray*}
0 &\lra& \chom(\cO_{\cA},\cI_{\cM}) \lra \chom(\cI_{\cM},\cI_{\cM}) \lra \chom(\cI_{\cZ},\cI_{\cM})\\
&\lra& \cext^1(\cO_{\cA},\cI_{\cM}) \lra \cext^1(\cI_{\cM},\cI_{\cM}) \lra \cext^1(\cI_{\cZ},\cI_{\cM})\\
&\lra& \cext^2(\cO_{\cA},\cI_{\cM}) \lra \cext^2(\cI_{\cM},\cI_{\cM}) \lra \cext^2(\cI_{\cZ},\cI_{\cM}) \lra 0.
\end{eqnarray*}

Note that the point $a$ is not in $M$. By Serre duality, it's easy to see that 
\begin{eqnarray*}
\cext^k(\cO_{\cA},\cI_{\cM}) &=& \cext^{2-k}(\cI_{\cM},\cO_{\cA}\otimes \omega_{\cX/\P^1})^\vee\\
&=&\cext^{2-k}(\cI_{\cM},\cO_{\cA})^\vee\otimes\cO_{\P^1}(2)\\
&=&\cext^{2-k}(\cO_{\cX},\cO_{\cA})^\vee\otimes\cO_{\P^1}(2)\\
&=&
\begin{cases}
\cO_{\P^1}(2) & \text{if $k=2$;}\\
0 & \text{otherwise.}
\end{cases}\\
\end{eqnarray*}

By the proof of Lemma \ref{case02} and Lemma \ref{samez} we have that 
$$\cext^k(\cI_{\cM},\cI_{\cM})=
\begin{cases}
\cO_{\P^1} & \text{if $k=0$;}\\
\cO_{\P^1}(1)^{\oplus (2p-2)} & \text{if $k=1$;}\\
\cO_{\P^1}(2) & \text{if $k=2$.}
\end{cases}
$$

Note that the natural map $$\chom(\cI_{\cM},\cI_{\cM}) \lra \chom(\cI_{\cM},\cO_{\cA})$$ given by the evaluation on $\cA$ is an isomorphism, therefore its dual map $$\cext^2(\cO_{\cA},\cI_{\cM}) \lra \cext^2(\cI_{\cM},\cI_{\cM})$$ is also an isomorphism. From the above long exact sequence we can conclude that 
$$\cext^k(\cI_{\cZ},\cI_{\cM})=
\begin{cases}
\cO_{\P^1} & \text{if $k=0$;}\\
\cO_{\P^1}(1)^{\oplus (2p-2)} & \text{if $k=1$;}\\
0 & \text{if $k=2$.}
\end{cases}$$

From another short exact sequence $$0 \lra \cI_{\cW} \lra \cI_{\cM} \lra \cO_{\cB} \lra 0$$ we can get
\begin{eqnarray*}
0 &\lra& \chom(\cI_{\cZ},\cI_{\cW}) \lra \chom(\cI_{\cZ},\cI_{\cM}) \lra \chom(\cI_{\cZ},\cO_{\cB})\\
&\lra& \cext^1(\cI_{\cZ},\cI_{\cW}) \lra \cext^1(\cI_{\cZ},\cI_{\cM}) \lra \cext^1(\cI_{\cZ},\cO_{\cB})\\
&\lra& \cext^2(\cI_{\cZ},\cI_{\cW}) \lra \cext^2(\cI_{\cZ},\cI_{\cM}) \lra \cext^2(\cI_{\cZ},\cO_{\cB}) \lra 0
\end{eqnarray*}

The point $b$ doesn't belong to $Z$, therefore we have
\begin{eqnarray*}
\cext^k(\cI_{\cZ},\cO_{\cB}) &=& \cext^k(\cO_{\cX},\cO_{\cB})\\
&=& \begin{cases}
\cO_{\P^1} & \text{if $k=0$;}\\
0 & \text{otherwise.}
\end{cases}
\end{eqnarray*}

Similarly as above, note that the natural map $$\chom(\cI_{\cZ},\cI_{\cM}) \lra \chom(\cI_{\cZ},\cO_{\cB})$$ given by the evaluation on $\cB$ is an isomorphism. Therefore from the above long exact sequence we can conclude that $$\cext^1(\cI_{\cZ},\cI_{\cW})=\cO_{\P^1}(1)^{\oplus (2p-2)}.$$

\end{proof}

Therefore we get the following corollary: 

\begin{corollary}\label{cor_zero}
Under the assumption of Lemma \ref{manysame}, we have
$$\ext^1_{\cX}(\cI_{\cZ},\cI_{\cW}\otimes\pi^*\cO_{\P^1}(-2))=0.$$
\end{corollary}

\begin{proof}
By local-to-global spectral sequence we have that 
\begin{eqnarray*}
0 &\lra& H^1(\chom(\cI_{\cZ},\cI_{\cW}\otimes\pi^*\cO_{\P^1}(-2))) \lra \ext^1_{\cX}(\cI_{\cZ},\cI_{\cW}\otimes\pi^*\cO_{\P^1}(-2))\\ 
&\lra& H^0(\cext^1(\cI_{\cZ},\cI_{\cW}\otimes\pi^*\cO_{\P^1}(-2))) \lra 0.
\end{eqnarray*}

First of all, on each K3 fiber of the twistor space, since $Z$ and $W$ are not identical, we have $$\hom(\cI_Z, \cI_W)=0,$$ hence $$\chom(\cI_{\cZ},\cI_{\cW}\otimes\pi^*\cO_{\P^1}(-2))=0.$$

Furthermore, by Lemma \ref{manysame},  $$\cext^1(\cI_{\cZ},\cI_{\cW}\otimes\pi^*\cO_{\P^1}(-2)) =\cO_{\P^1}(1)^{\oplus (2p-2)}\otimes\cO_{\P^1}(-2)=\cO_{\P^1}(-1)^{\oplus(2p-2)},$$ so we also have $$H^0(\cext^1(\cI_{\cZ},\cI_{\cW}\otimes\pi^*\cO_{\P^1}(-2)))=0.$$ Therefore, we conclude that $$\ext^1_{\cX}(\cI_{\cZ},\cI_{\cW}\otimes\pi^*\cO_{\P^1}(-2))=0.$$
\end{proof}

Before moving on to the third step, we need to collect some known results about sections of the twistor fibration $$\pi:\cX \lra \P^1.$$ Note that every section of this fibration is isomorphic to $\P^1$. The space of sections of the twistor fibration is a complex variety. More restrictively, if we only consider sections whose normal bundle is $\cO_{\P^1}(1)^{\oplus 2}$, then we get a smooth open subvariety, which we denote by $\sec$. 

All the horizontal sections of the twistor fibration form a real manifold $\hor$, which is the underlying real manifold of the K3 surface we start with. Since every horizontal section has a normal bundle $\cO_{\P^1}(1)^{\oplus 2}$, $\hor$ is a real analytic submanifold of $\sec$. However, it doesn't inherit a complex structure from $\sec$. 

Since both $\cZ$ and $\cW$ are unions of $p$ horizontal sections, we can think of them as closed points in $\hor^p$ or $\sec^p$. We denote by $\sec^p_0$ the open subvariety consisting of $p$ disjoint sections, each of which has a normal bundle $\cO_{\P^1}(1)^{\oplus 2}$. Similarly, we denote by $\hor^p_0$ the open submanifold consisting of $p$ disjoint horizontal sections. It's easy to see that $\hor^p_0$ is exactly the underlying real analytic manifold of $X^p_0$. 

Now we are ready to finish the proof of Proposition \ref{diffz}. 

\begin{proof}[Proof of Proposition \ref{diffz}]

Let $\mathcal{C}$ be the universal family of $p$ disjoint sections over $\sec^p_0$. Consider the product $\sec^p_0\times\sec^p_0$ with projections to two factors $\pi_1$ and $\pi_2$. Then $\pi_1^*\mathcal{C}$ and $\pi_2^*\mathcal{C}$ are both flat over $\sec^p_0\times\sec^p_0$. 

By a version of upper semi-continuity theorem \cite[Theorem 3]{Ban-german}, the dimension of $$\ext^1_{\cX}(\cI_{\cZ},\cI_{\cW}\otimes\pi^*\cO_{\P^1}(-2))$$ is upper semi-continuous on $\sec^p_0\times\sec^p_0$, where $(\cZ,\cW)\in\sec^p_0\times\sec^p_0$. 

However, by Corollary \ref{cor_zero}, we know that at some closed points of $\hor^p_0\times\hor^p_0$, we have $$\ext^1_{\cX}(\cI_{\cZ},\cI_{\cW}\otimes\pi^*\cO_{\P^1}(-2))=0.$$ Therefore, by the upper semi-continuity, there is an open subset of $\sec^p_0\times\sec^p_0$, on which the same equation holds. We denote the intersection of this open subset and $\hor^p_0\times\hor^p_0$ by $V$, then $V$ is a non-empty analytic open subset of $\hor^p_0\times \hor^p_0$. For any $(\cZ,\cW)\in V$, we look at the local-to-global spectral sequence
\begin{eqnarray*}
0 &\lra& H^1(\chom(\cI_{\cZ},\cI_{\cW}\otimes\pi^*\cO_{\P^1}(-2))) \lra \ext^1_{\cX}(\cI_{\cZ},\cI_{\cW}\otimes\pi^*\cO_{\P^1}(-2))\\ 
&\lra& H^0(\cext^1(\cI_{\cZ},\cI_{\cW}\otimes\pi^*\cO_{\P^1}(-2))) \lra 0.
\end{eqnarray*}
Since the middle term is 0, we conclude that the right term is also 0. 
However, $$\cext^1(\cI_{\cZ},\cI_{\cW}\otimes\pi^*\cO_{\P^1}(-2)) =\cext^1(\cI_{\cZ},\cI_{\cW})\otimes\cO_{\P^1}(-2).$$
Therefore, we conclude that $\cext^1(\cI_{\cZ},\cI_{\cW})$ does not have any direct summand $\cO_{\P^1}(2)$. By Lemma \ref{lemma_ext}, we finally get $$\cext^1(\cI_{\cZ},\cI_{\cW})=\cO_{\P^1}(1)^{\oplus (2p-2)},$$ where $\cZ$ and $\cW$ are the twistor lift of $Z$ and $W$, which consist of $p$ distinct points each. 

\end{proof}

Now we can finish the proof of Proposition \ref{mainprop}. 

\begin{proof}[Proof of Proposition \ref{mainprop} (2)]

By the argument at the beginning of this section, it suffices to prove the formality of $\rhom(E,E)$ for $$E=\bigoplus_{i=1}^m\cI_{Z_i}^{n_i},$$ where each $Z_i\in X^p_0$ is the union of $p$ distinct points in $X$, and all $Z_i$'s are distinct. We denote the big diagonal in the product $(X^p_0)^m$ by $\Delta$, then $$(Z_1, Z_2, \cdots, Z_m)\in (X^p_0)^m\backslash\Delta.$$

By Proposition \ref{diffz}, we can find an open subset $$U\subset (X^p_0)^m\backslash\Delta$$ in classical topology, such that for any $$(Z_1, Z_2, \cdots, Z_m)\in U,$$ we have $\cext^1(\cI_{\cZ_i}, \cI_{\cZ_j})$ has pure weight 1. Together with Proposition \ref{samez}, we see that $\cext^1(\cE, \cE)$ has pure weight 1. Together with Lemma \ref{case02}, we see that for every $k$, $\cext^k(\cE, \cE)$ has pure weight $k$. By the same argument as in the proof of Proposition \ref{mainprop} (1), we apply Theorem \ref{kalthm} (ii) to conclude that $\rhom(E,E)$ is formal in this case. 

For a general choice of $$(Z_1, Z_2, \cdots, Z_m)\in (X^p_0)^m\backslash\Delta$$ which is not necessarily in $U$, we need to apply Theorem \ref{kalthm} (i). Note that $U$ is an open subset of $(X^p_0)^m\backslash\Delta$ in classical topology, not in Zarisky topology. However, the proof of Theorem \ref{kalthm} (i) in \cite{Kal-family} can be carried over in this case without any change. Therefore, since the formality holds in $U$, it also holds in $(X^p_0)^m\backslash\Delta$, which is what we want to prove. 

\end{proof}

\section{The Remaining Case}

From the above discussion, we see that the most important point in our proof is that the relative extension sheaf $\cext^k(\cE, \cE)$ is of pure weight $k$. In the previous two sections we basically showed that, if all direct summands in the decomposition of $E$ are simultaneously hyperholomorphic bundles or simultaneously ideal sheaves of 0-dimensional subschemes, the relative extension sheaf $\cext^k(\cE, \cE)$ has the expected weight. However, if both hyperholomorphic bundles and ideal sheaves appear in the stable decomposition of $E$, the following lemma shows that $\cext^1(\cE, \cE)$ doesn't have pure weight 1 any more. In fact, it acquires summands of $\cO_{\P^1}$, $\cO(1)_{\P^1}$ and $\cO(2)_{\P^1}$. 

\begin{lemma}
Let $\cF$ and $\cI_{\cZ}$ be twistor transforms of a stable hyperholomorphic bundle and an ideal sheaf of $p$ distinct points respectively. Then $\cext^1(\cF,\cI_{\cZ})$ acquires direct summands of both $\cO_{\P^1}$ and $\cO_{\P^1}(1)$, and $\cext^1(\cI_{\cZ},\cF)$ acquires direct summands of both $\cO_{\P^1}(1)$ and $\cO_{\P^1}(2)$.
\end{lemma}

\begin{proof}
By resolving the ideal sheaf of horizontal sections $\cI_{\cZ}$, we have the long exact sequence
\begin{eqnarray*}
0 &\lra& \chom(\cF,\cI_{\cZ}) \lra \chom(\cF,\cO_{\cX}) \lra \chom(\cF,\cO_{\cZ})\\
&\lra& \cext^1(\cF,\cI_{\cZ}) \lra \cext^1(\cF,\cO_{\cX}) \lra \cext^1(\cF,\cO_{\cZ})
\end{eqnarray*}
Without loss of generality, we can additionally assume that $\cF$ is $\mu$-stable, then $$\chom(\cF,\cO_{\cX})=0.$$ In fact, its restriction on each fiber vanishes, since we have assumed $\rk F$ has rank at least 2. By the argument of Lemma \ref{dualhyper}, we know that $F^\vee$ is also a hyperholomorphic bundle. Using Theorem \ref{twistlift}, we get 
$$\chom(\cF,\cO_{\cZ}) =\chom(\cO_{\cX},\cF^\vee_{|\cZ}) =\pi_*\cF^\vee_{|\cZ}=\pi_*\cO_{\cZ}^{\oplus r}=\cO_{\P^1}^{\oplus rp}$$
where $r$ is the rank of $F$ and $p$ is the length of $Z$. 
Moreover, $$\cext^1(\cF,\cO_{\cX}) =\cext^1(\cO_{\cX},\cF^\vee) =R^1\pi_*\cF^\vee=\cO_{\P^1}(1)\otimes H^1(F^\vee).$$ Note that in the last step we used Proposition \ref{weight}. Finally, we have $$\cext^1(\cF,\cO_{\cZ}) =\cext^1(\cO_{\cX},\cF^\vee_{|\cZ}) =R^1\pi_*\cF^\vee_{|\cZ}=R^1\pi_*\cO_{\cZ}^{\oplus r}=0.$$
Therefore the exact sequence becomes
$$0 \lra \cO_{\P^1}^{\oplus rp} \lra \cext^1(\cF,\cI_{\cZ}) \lra \cO_{\P^1}(1)\otimes H^1(F^\vee) \lra 0.$$ Since $\ext^1_{\P^1}(\cO_{\P^1}(1),\cO_{\P^1})=0$, we know that $$\cext^1(\cF,\cI_{\cZ})=\cO_{\P^1}^{\oplus rp}\oplus \cO_{\P^1}(1)\otimes H^1(F^\vee)$$ which has a mixed weight of 0 and 1.

Similarly, it's not hard to see that $$ 0 \lra \cext^1(\cO_{\cX},\cF) \lra \cext^1(\cI_{\cZ},\cF) \lra \cext^2(\cO_{\cZ},\cF) \lra 0,$$ where $\cext^1(\cO_{\cX},\cF)$ has pure weight 1, by Proposition \ref{weight}, and $\cext^2(\cO_{\cZ},\cF)$ has pure weight 2, by Theorem \ref{twistlift} and relative Serre duality. Therefore $\cext^1(\cI_{\cZ},\cF)$ is the direct sum of the two, which has a mixed weight 1 and 2.
\end{proof}

Therefore, in this mixed case, Theorem \ref{kalthm} (ii) doesn't apply anymore. It seems that we need new techniques to deal with it. I hope to come back to this point later.

\bibliographystyle{alpha}
\bibliography{references}

\end{document}